\input epsf
\documentclass{amsart}
\usepackage{amssymb,latexsym,hyperref,epsfig}

\numberwithin{equation}{section}
\newtheorem{theorem}{Theorem}[section]

\newtheorem{remark}[theorem]{Remark}
\newtheorem{lemma}[theorem]{Lemma}

\newtheorem{example}[theorem]{Example}

\newtheorem{claim}[theorem]{Claim}

\def\mM{{\mathcal M}}

\def\pP{{\mathcal P}}

\def\ZZ{{\mathbb Z}}

\newenvironment{proofof}[1]{\begin{trivlist}\item {\emph{Proof of {#1}.}\,}}
{\mbox{}\hfill$\square$\end{trivlist}}

\begin{document}

\title
{Existence of perfect Morse functions on
spaces with semi-free circle action.}

\author {Mikhail Kogan}

\thanks{The author is supported by NSF Postdoctoral Fellowship.}

\address{\noindent Department of Mathematics, Northeastern University,
  Boston, MA 02115}
\email{misha@neu.edu}

\begin{abstract} Let $M$ be a compact oriented simply-connected manifold of
dimension~at least $8$. Assume $M$ is equipped with a torsion-free
semi-free circle action with isolated fixed points. We prove $M$
has a perfect invariant Morse-Smale function. The major ingredient
in the proof is a new cancellation theorem for the invariant Morse
theory.
\end{abstract}

\maketitle

\pagestyle{myheadings} \markboth{MIKHAIL KOGAN}{PERFECT MORSE
FUNCTIONS ON SPACES WITH SEMI-FREE CIRCLE ACTIONS}

\section{Introduction.}

It is well known that the classical results of Morse theory
\cite{m,s} do not directly generalize to the invariant case. In
particular, a compact manifold with a compact group action does
not necessarily have a Morse-Smale function, that is a Morse
function for which stable and unstable manifolds intersect
transversely. As a result, in the presence of a group action there
is no guarantee of existence of perfect invariant Morse function
for which Morse inequalities become equalities making the
invariant Morse polynomial of the function equal to the invariant
Poincare polynomial of the manifold.

At the same time, there are plenty of examples, many of which
appear in symplectic geometry, of manifolds with groups actions
which carry perfect invariant Morse functions. As shown by Atiyah
and Bott \cite{ab,ab1} all compact symplectic manifolds with
Hamiltonian torus action are equipped with perfect invariant Morse
functions which arise as components of moment maps.

The existence of perfect Morse function is one of the many
topological properties of Hamiltonian spaces for which it is
natural to ask about generalization to non-Hamiltonian actions.
Some of these properties, the Kirwan surjectivity theorem
\cite{ki} and Jeffrey-Kirwan localization formula \cite{jk}, have
been studied by the author \cite{k}. It was shown that
Jeffrey-Kirwan localization formula can be extended to the case of
equivariantly formal torus actions, while Kirwan surjectivity
theorem does not generalize to this case.

This paper, an attempt to generalize the existence of perfect
invariant Morse functions to non-Hamiltonian actions, was
motivated by the recent result of Hattori \cite{h} and
Tolman-Weitsman \cite{tw} on semi-free circle actions on
symplectic manifolds. Tolman and Weitsman showed that every
semi-free symplectic circle action on a compact symplectic
manifold with finite nonempty fixed point set must be Hamiltonian
and hence has a moment map, a perfect invariant Morse function.

Our main theorem generalizes the result to spaces with semi-free
action whose equivariant cohomology is $\ZZ$-torsion-fee. For the
definitions of all the terms used in the statement of the theorem
see Section \ref{definition}.




\begin{theorem}
\label{perfect-morse} Let $M$ be a compact oriented
simply-connected manifold with \mbox{torsion-free}
\mbox{semi-free} circle action which has finitely many fixed
points. Further assume \mbox{$\dim M\geq 8$}. Then there exists a
perfect invariant Morse-Smale function on $M$.
\end{theorem}

While the original motivation for this result comes from the
symplectic geometry, the methods of proof are borrowed from the
classical Morse theory, namely the ideas used in the proof of
$h$-cobordism theorem in \cite{m}. The major ingredient in the
proof is the new cancellation theorem
(Theorem~\ref{cancell-theorem}) which is similar to the The
cancellation theorem supplies the tools for eliminating extra
critical circles of invariant Morse functions.

\begin{figure}[ht]
\epsffile{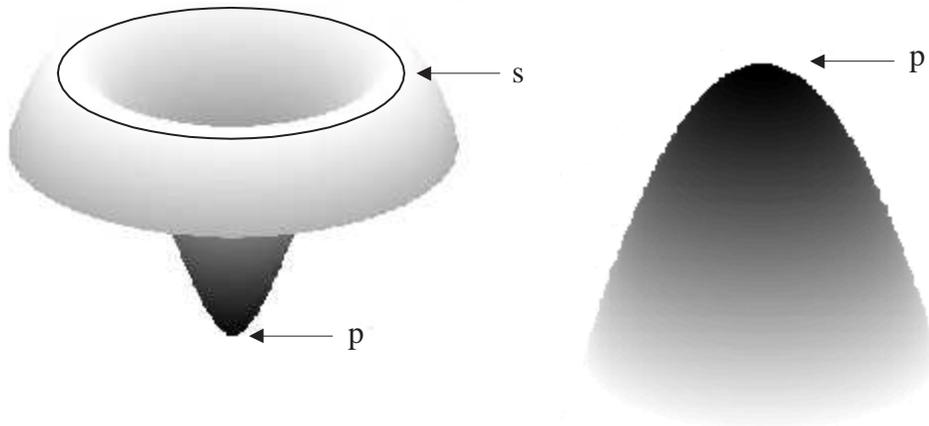} \label{figure} \caption{Critical circle $s$
is cancelled by increasing the index of critical point $p$ from
$0$ to $2$.}
\end{figure}

The essential idea behind the cancellation theorem is illustrated
Figure~\ref{figure}, where a critical circle of index $1$ is
cancelled by increasing the index of a critical point from $0$
to~$2$. The action of the circle on both pictures in
Figure~\ref{figure} is given by rotation around the vertical axis,
while the Morse function is given by projection onto the same
axis. In the first picture there are a critical point $p$ of index
$0$ and a critical circle $s$ of index~$1$, while in the second
picture $s$ is eliminated by increasing the index of $p$ to $2$.

\vspace{-2mm}

The paper is organized as follows. In Section \ref{definition} we
provide basic definitions and notations, and state the main
results of the paper. Section \ref{section-morse} contains the
proof of existence of Morse-Smale functions on spaces with
semi-free circle actions (Theorem~\ref{morse-smale}). Section
\ref{section-cancell} provides the proof of cancellation theorem
(Theorem~\ref{cancell-theorem}). The connections between invariant
Morse theory and equivariant cohomology are discussed in Section
\ref{cohomology}. Finally, the last section is concerned with the
proof of Theorem \ref{perfect-morse}, which boils down to
eliminating critical circles using the cancellation theorem.

Let us make a comment about the presentation. Since Morse theory
is a vast subject, we can not mention all contributors and refer
to all known results. So, to make our presentation comprehensible,
we quote mainly one source (the book of Milnor~\cite{m}).
To~shorten the presentation we often refer to certain statements
or proofs in~\cite{m} and indicate the how to change them to
produce analogous statements and proofs in the presence of the
circle action.

\section{Definitions and Main Results.}
\label{definition}

A smooth action of a circle $S^1$ on a manifold $M$ is called {\it
semi-free} if the action is free outside the set of fixed points
$M^{S^1}$. In this paper we always assume that the action of~$S^1$
on~$M$ is semi-free and there are finitely many fixed points of
the action.

An important property of semi-free actions is that around every
isolated fixed point~$p$, there exists an open neighborhood~$U$ of
$p$ and an invariant diffeomorphism $\phi$ between~$U$ and an
invariant open subset of $\mathbb C^n$, where $S^1$ acts on
$\mathbb C^n$ with weight~$(1,\dots,1)$. We call~$(U,\phi)$ a {\it
standard chart} at $p$. We will use both complex coordinates
$(z_1,\dots, z_{2n})$ and real coordinates
$(x_1,y_1,\dots,x_{2n},y_{2n})$ on $\mathbb C^n=\mathbb R^{2n}$
with $z_i=x_i+\sqrt{-1} y_i$.

Let $f:M^{2n}\to \mathbb R$ be a smooth invariant function. Denote by
$Crit(f)$ the set of critical points of $f$. Let $Crit_c(f)$
be the set of connected components of
$Crit(f)$. Notice that all isolated fixed point of the action must be
critical points of $f$.

A fixed point $p\in M^{S^1}$ is {\it nondegenerate} if the
Hessian of $f$ at $p$ is nondegenerate. An invariant version of
the Morse lemma says that for every nondegenerate fixed point
there exists a standard chart $(U,\phi)$ at $p$, so that on $U$
$$
f=\phi^*(f(p)-|z_1|^2-\dots-|z_\lambda|^2
+|z_{\lambda+1}|^2+\dots+|z_n|^2).
$$
The chart $(U,\phi)$ is called {\it a standard chart of $f$ at
$p$}, the index $\sigma(p)$ of $p$ is set to be~$2\lambda$.

A circle orbit $s\subset Crit(f)$ will be called \emph{a critical
circle of $f$}. Choose an
invariant
neighborhood $\tilde U$ of $s$ on which $S^1$ acts freely and let
$\pi$ be the map $\tilde U\to \tilde U/S^1$. Then $s$ is {\it
nondegenerate} if the Hessian of~$\pi_*(f)$ is nondegenerate at
$s/S^1$. The Morse lemma states that there exists an invariant
neighborhood $U\subset \tilde U$ of $s$ and coordinates
$(x_1,...,x_{2n-1})$ on $U/S^1$, such that
$$
f=\pi^*(f(p)-x_1^2-\dots-x_\lambda^2
+x_{\lambda+1}^2+\dots+x_{2n-1}^2).
$$
We say that $U$ together with the coordinates $x_1,\dots,x_{2n-1}$
is {\it a standard chart of $f$ at $s$}. The index $\sigma(s)$ of
$s$ is defined to be $\lambda$.

Define $f$ to be {\it an invariant Morse function} if each connected
component of $Crit(f)$ is either a nondegenerate fixed point or
a nondegenerate critical circle.

Assume $f$ is an invariant Morse function on $M$. An invariant
vector field $\xi$ on $M$ is {\it a gradient-like vector field} if
\begin{enumerate}
\item{$\xi(f)<0$ on $M-Crit(f)$,}
\item{for every nondegenerate critical point $p$, there exists a
standard chart $(U,\phi)$ of~$f$ at~$p$, such that $\xi$ on $U$
has coordinates
$$
(x_1,y_1,\dots,x_\lambda,y_\lambda,-x_{\lambda+1},-y_{\lambda+1},
\dots,-x_n,-y_n).
$$}
\item{for every nondegenerate critical circle $s$ there exists a
standard chart of $f$ at~$s$, such that the push forward vector
field $\pi_*(\xi)$ on $U/S^1$ has coordinates
$$
(x_1,\dots,x_\lambda,-x_{\lambda+1},\dots,-x_{2n-1}).
$$}
\end{enumerate}


A gradient-like vector field $\xi$ defines a smooth flow
$\phi_\xi(t)$ on $M$. For \mbox{$q\in Crit_c(f)$}~define
$$
W^\pm_q=\{x\in M|\lim_{t\to\pm\infty}\phi_\xi(t)\in q\}
$$
The manifold $W^+_q$ is the {\it stable manifold} of $q$, and
$W^-_q$ is the {\it unstable manifolds} of $q$. The tuple
$(f,\xi)$ is called {\it an invariant Morse-Smale function} on $M$
if every stable manifold intersects every unstable manifolds
transversely.

\begin{theorem}
\label{morse-smale} Every semi-free circle action on a compact
manifold with isolated fixed points has an invariant Morse-Smale
function.
\end{theorem}

The key step in the proof of our main result is a new cancellation
theorem. To state it, let us introduce some notations. For a
subset $I$ of $\mathbb R$ and a submanifold $N$ of $M$, define
$N_I(f)=f^{-1}(I)\cap N$. Whenever no confusion may arise, we will
drop $f$ from the notation and set $N_I=N_I(f)$. For example, for
$a\in \mathbb R$, $W^+_{q,a}=W^+_q\cap f^{-1}(a)$, $W^-_{q,\leq
a}=W^-_q\cap f^{-1}((-\infty, a])$, or for $a<b$,
$M_{[a,b]}=f^{-1}([a,b])$.

Assume $V$ is an oriented manifold with a free $S^1$ action,
$N,N'$ are two submanifolds of $V$ which intersect transversely
with $\dim N+\dim N'=\dim V+1$, and $N$ and the normal bundle $\nu
(N')$ are oriented. In Section \ref{section-cancell} we define
{\it the intersection number $N\cdot N'$} by counting the number
of circles in $N\cap N'$ with appropriate signs. Orientations
on~$N$ and $\nu(N)$ are {\it compatible} with respect to an
orientation on $V$, if at any point of $N$ a positively oriented
frame on $N$ followed by a positively oriented frame on $\nu(N)$
gives a positively oriented frame on $V$. If the orientation on
$N$ and $\nu(N)$ as well as on $N'$ and~$\nu(N')$ are compatible
to an orientation on $V$, then $N\cdot N'=\pm N'\cdot N$.

\begin{theorem}
\label{cancell-theorem} {\rm (The cancellation theorem.)} Assume
the action of $S^1$ on a compact manifold $M$ is semi-free and
$(f,\xi)$ is an invariant Morse-Smale function on $M$. Assume~$a<b$
and $M_{(a,b)}$ contains one fixed point $p$ of index
$\lambda$ and one critical circle $s$ of index~$\lambda+1$. Let
$f(p)<c<f(s)$. Suppose $V=M_c$ is oriented and $N=W^+_{p,c}$,
$\nu(N)$ as well as $N'=W^-_{s,c}$, $\nu(N')$ are equipped with
compatible orientations and $N \cdot N'=\pm 1$. Further assume
that either $\dim M\geq 8$, $M_a/S^1$ is simply-connected, and
$M_a$ is connected; or $\lambda=0$; or $\lambda=\dim M-2$.

Then there exist a new invariant Morse function $f'$ and a new
gradient-like vector field $\xi'$ identical to $(f,\xi)$ outside
$M_{(a,b)}$ such that $M_{(a,b)}$ contains no critical circles
of $f'$ and the index of the fixed point $p$ is $\lambda+2$.
\end{theorem}

A circle action is \emph{torsion-free} if the equivariant
cohomology \mbox{$H^*_{S^1}(M)=H^*_{S^1}(M, \ZZ)$} is torsion-free
as a $\ZZ$ module. For a torsion-free $S^1$ action define its
\emph{Poincare polynomial}~by
\begin{eqnarray*}
\pP_{S^1}(M) & = & \sum_{k\in \ZZ_{\geq 0}} \dim
(H^k_{S^1}(M))t^k.
\end{eqnarray*}
For a fixed point point $p$ define
$\mM(p)=\frac{t^{\sigma(p)}}{1-t^2}$, and for a critical
circle~$s$ set $\mM(s)=t^{\sigma(s)}$. Define \emph{the Morse
polynomials} of an invariant Morse function $f$ by
\begin{eqnarray*}
\mM_{S^1}(M,f) & = & \sum_{q\in Crit_c(f)} \mM(q).
\end{eqnarray*}
An invariant Morse function is {\it perfect} if the Poincare and
Morse polynomials are the same $\pP_{S^1}(M)=\mM(M,f)$.
We have given all the definitions needed to state
Theorem~\ref{perfect-morse}.

\section{Existence of invariant Morse-Smale functions}
\label{section-morse}

We believe many results of this section are known to the experts
in the field (see~\cite{b} for further references). We provide the
proofs of these results for completeness.

\begin{theorem}
\label{morse} Every semi-free circle action on a compact manifold
with isolated fixed points has an invariant Morse function.
\end{theorem}

\proof Let $M^{S^1}=\{p_1,\dots,p_m\}$. Around every fixed point
$p_i$ there exists a standard coordinate chart $(U_i,\phi_i)$.
Without loss of generality, we can assume each $U_i$ is
diffeomorphic to a unit ball $B_1$ inside of $\mathbb C^n$ and
that $U_i$'s are pairwise disjoint.

Define $f=\phi_i^*(|z|^2)$ on every $U_i$. We want to extend $f$ to
$M-U$ , where $U=\cup_iU_i$. The manifold $X=(M-U)/S^1$ is compact
and $\partial X=(\partial \bar U) /S^1$. By \cite[Theorem~2.5]{m}
there exists a Morse function $h:(X,\partial X)\to ([1,\infty),1)$.
Denote by $\pi$ the projection~\mbox{$M- U\to X$}. Define
$f=\pi^*(h)$ on
$M-U$. Then $f$ is smooth outside of $\partial \overline U$.

To finish the proof use the following smoothing argument. Choose
an invariant tubular neighborhood $V$ of $\partial \bar U$ and
identify it with $(-\varepsilon, \varepsilon)\times \partial \bar
U$. Let $f^-$, $f^+$ be two smooth functions on $V$, such that
$f^-(t,x)=f(t,x)$ for $(t,x)\in (-\varepsilon, 0)\times \partial
\bar U$ and $f^+(t,x)=f(t,x)$ for $(t,x)\in (0,\varepsilon)\times
\partial \bar U$. Without loss of generality, we can assume
$\frac{\partial f^-}{\partial t}(t,x)>0$ and $\frac{\partial
f^+}{\partial t}(t,x)>0$ on $V$ and $f^+(t,x)>f^-(t,x)$ for
$|t|>\varepsilon/2$. Let $\mu(t)$ be a smooth positive
nondecreasing function, with $\mu(t)=0$, for $t\leq -\varepsilon$;
$ \mu(t)=1$ for $t\geq \varepsilon$; and $\mu(t)=\frac{1}{2}$ for
$|t|<\varepsilon/2$. Change $f$ to $\mu f^+ +(1-\mu) f^-$ on $V$.
Then $f$ is smooth on $M$. An easy computation shows that
$\frac{\partial f}{\partial t} >0$ on $V$ which implies that this
construction does not create new critical points of $f$.
\endproof

Every invariant Morse function together with an invariant
gradient-like vector field on a manifold with isolated fixed
points defines a partial order on $Crit_c(f)$ as follows. If
$q_1,q_2\in Crit_c(f)$, set $q_1\prec q_2$ if $q_1$ is in the
closure of $W_{q_2}^{-}$. Complete this relation to satisfy
transitivity and hence define a partial order on $Crit_c(f)$. A
map $c:Crit_c(f)\to \mathbb R$ is said to be {\it order
preserving} if $q_1\prec q_2$ implies $c(q_1)<c(q_2)$.

\begin{theorem}
\label{order} Assume the circle action on $M$ is semi-free and has
isolated fixed points. Let $f$ be an invariant Morse function and
$\xi$ be its gradient-like vector field on $M$, then every order
preserving map on $Crit_c(f)$ can be extended to an invariant
Morse function.
\end{theorem}

\proof Throughout the proof $\xi$ and the partial order $\prec$
are fixed. Notice that for every~\mbox{$q\in Crit_c(f)$} there
exists an invariant Morse function~$g$, such that $g=f$ outside a
small invariant neighborhood  $U$ of $q$, $g(q)=f(q)+\varepsilon$
for a small enough $\varepsilon$, and $\xi$ is an invariant
gradient-like vector field for $g$. The function $g$ can be
defined by $g=f+\rho$, where $\rho$ is an invariant ``bump''
function, such that $\rho=0$ outside of $U$, $\rho=\varepsilon$ on
a compact set~$K\subset U$ and $|\xi(\rho)|<|\xi(f)|$ on $U-K$.

Hence the theorem follows from the following statement. Let~$f$
have two critical values in the interval $(a,b)$,
such that one critical level contains \mbox{$\{q_1,\dots,q_k\}\subset
Crit_c(f)$} and the other contains $\{r_1,\dots,r_m\}\subset
Crit_c(f)$. Moreover, assume non of the~$q_i$'s are compatible
with $r_j$'s in the partial order $\prec$, that is stable
and unstable manifolds of $q_i$'s and $r_j$'s never intersect.
Then for any $a_1, a_2 \in (a,b)$, there exists a invariant
Morse function~$g$, such that $\xi$ is an invariant gradient-like
vector field of $g$, $f=g$ outside $M_{[a,b]}$, the set of
critical points of $g$ is the same as the set of critical points
of $f$, and $g(q_i)=a_1$, $g(r_j)=a_2$.

Without the word ``invariant'' the above statement is
\cite[Theorem~4.2]{m}. This allows to give a proof of this statement
by making a few minor changes in the proof of
\mbox{\cite[Theorem~4.2]{m}},
mainly by inserting the word ``invariant'' where necessary.
\endproof

Let $S^1$ act on two manifolds $M$ and $N$. \emph{An invariant
isotopy} between two invariant diffeomorphisms $h_0,h_1:M\to N$ is
a smooth map $h:M\times[0,1]\to N$, such that $h_0(x)=h(x,0)$,
$h_1(x)=h(x,1)$ and every $h_t(x)=h(x,t)$ is an invariant
diffeomorphism. We say that $h$ is \emph{an isotopy of identity}
if $M=N$ and $h_0$ is an identity map.

\begin{lemma}
\label{skew} Assume we are given an invariant Morse function $f$
on $M$ with an invariant gradient-like vector field $\xi$, a
non-critical level $M_a=f^{-1}(a)$, and a diffeomorphism
$h_1:M_a\to M_a$ that is invariantly isotopic to the identity. If
for $a<b$ the interval $[a,b]$ contains no critical values, then
it is possible to construct a gradient-like vector field $\bar
\xi$, such that $\xi=\bar \xi$ outside $M_{[a,b]}$ and
$\bar\varphi=h_1\circ\varphi$, where $\varphi$ and $\bar \varphi$
are the diffeomorphisms $M_b\to M_a$ determined by following the
trajectories of $\xi$ and $\bar \xi$ respectively.
\end{lemma}

\proof This lemma is the invariant version of \cite[Lemma~4.7]{m}
whose proof can be modified (mostly by inserting the word
``invariant'') to provide a proof of the lemma.
%
%
\endproof

\begin{proofof}{Theorem \ref{morse-smale}}
Using Theorem~\ref{morse} choose an invariant Morse function $f$
on~$M$. Pick any invariant riemmannian metric on $M$, let
$\nu=-grad (f)$. On canonical charts~$U_i$ of $f$ at every fixed
point or critical circle define a vector field $\mu$ which is a
gradient-like vector field on $\cup U_i$. Pick a nonnegative
function $\rho$ which is zero outside of $\cup U_i$ and one in a
neighborhood of $Crit(f)$. Then $\xi=\rho\mu+(1-\rho )\nu$ is a
gradient-like vector field.

By Theorem~\ref{order}, we can assume that every critical level of
$f$ contains exactly one fixed point or one critical circle.
Let $c_1<\dots<c_k$ be the critical values of~$f$ and $q_1,\dots,
q_k \in Crit_c(f)$ with $f(q_i)=c_i$. By induction on $i$ we will
prove that there exists an invariant gradient-like vector field
$\xi$ of $f$, such that every $W^-_{q_j}$ intersects transversely
with every stable manifold for any $j\leq i$.

The base of the induction is trivial for $i=1$. Assume our
induction statement holds for $i-1$. To prove the theorem it is
enough to perturb $\xi$ on the preimage of $(c_{i-1},c_i)$ to
guarantee that $W^-_{q_i}$ intersects transversely with all the
stable manifolds.

Choose a number $a$ between $c_{i-1}$ and $c_i$. There exists an
invariant tubular neighborhood $U$ of~$W^-_{q_i,a}$ inside $M_a$
which is diffeomorphic to an $\mathbb R^k$ vector bundle $\pi:V\to
W^-_{q_i,a}$. If $q_i$ is a critical circle, then
$V$ can be chosen to be an $S^1$-trivial bundle, that is~$U$ is
invariantly diffeomorphic to $W^-_{q_i,a}\times \mathbb R^k$, with
$S^1$ acting only on the first component. Let $P:U\to \mathbb R^k$
be the natural projection. Let $Y_j$ be the set of critical values
of~$P$ restricted to $W^+_{q_j,a}\cap U$ and set $Y=\cup_jY_j$. By
Sard's theorem $Y$ has measure zero inside $\mathbb R^k$. Choose
$v\in \mathbb R^k-Y$. There exists an invariant isotopy of
identity $h_t$ on $U\cong W^-_{q_i,a}\times \mathbb R^k$, such
that $h_1(W^-_{q_i,a}\times 0)=W^-_{q_i,a}\times v$ and $h_t$ is
the identity outside some compact neighborhood $K\subset U$ of
$W^-_{q_i,a}$. Extend this isotopy to $M_a$ by setting it to be
identity outside of $U$. Then apply Lemma \ref{skew} to perturb
the gradient-vector field~$\xi$ to guarantee that the induction
assumption is satisfied.

If $q_i$ is a fixed point, then $\pi:V\to
W^-_{q_i,a}$ is no longer $S^1$-trivial. But~$W_{q_i,a}$ can be
covered by finitely many compact invariant sets
$\{K_\ell\}_{\ell=1}^L$ such that each~$K_\ell$ is inside an open
invariant set $W_\ell$ which invariantly contracts to a circle.
Without loss of generality we may assume that $W_\ell$'s are
invariantly diffeomorphic to $B_1^\circ\times S^1$ (where
$B^\circ_1$ is an open ball of radius 1 and $S^1$ acts only on the
second factor) and $K_\ell$'s are $B_{\frac{1}{2}}\times S^1$
inside $W_\ell$. Set $V_\ell=\pi^{-1}(W_\ell)$, then the bundle
$V_\ell\to W_\ell$ can be trivializes and, in particular,
each~$V_\ell$ is invariantly diffeomorphic to $W_\ell\times \mathbb
R^k$ with $S^1$ acting only on the first factor.

Let $P_1$ be the natural projection from $V_1$ to $\mathbb R^k$.
Assume $v_1$ is a regular value for every restriction of $P_1$ to
$W^+_{q_j,a}$. Then construct an isotopy of identity $h^1_t$ on
$M_a$ which is identity outside $V_1$ and such that $h_1^1$ takes
$K_1\times 0$ onto $K_1\times v_1$. This will guarantee that
$h^1_1(W^-_{q_i,a})$ is transverse to every $W^+_{q_j,a}$ on
$K_1$.

To construct an isotopy of identity $h^2_t$, choose a
trivialization \mbox{$h^1_1(V_2)=h_1^1(W_2)\times \mathbb R^k$},
let $P_2$ be the natural projection onto the second factor, and
choose $v_2$ which is a regular value for every restriction of
$P_2$ to $W^+_{q_j,a}$. Choose $h^2_t$, such that it is identity
outside of $h^1_1(V_2)$ and $h^2_1$ takes $h^1_1(K_2)\times 0$
onto $h^1_1(K_2)\times v_2$. Since~$K_1$ is compact and
transversality is an open condition we can choose~$v_2$ to be
small enough so that $h^2_1(h_1^1(W^-_{q_i,a}))$ is transverse to
$W^+_{q_j,a}$ on~$K_1$ Then $h^2_1(h_1^1(W^-_{q_i,a}))$ is
transverse to all~$W^+_{q_j,a}$ on $K_1\cup K_2$.

This process can be continued to construct invariant isotopies of
identity $h^m_t$, such that
$h^m_1(\dots(h_1^1(W^-_{q_i,a}))\dots)$ is transverse to every
$W^+_{q_j,a}$ on $\cup_{\ell=1}^m K_\ell$. Eventually we will
construct $h_t^1,\dots, h^L_t$ such that
$h^L_1(\dots(h_1^1(W^-_{q_i,a}))\dots)$ is transverse to every
$W^+_{q_j,a}$. Repeatedly apply Lemma \ref{skew} to perturb $\xi$
to a new gradient-like vector field to guarantee~$W^-_{q_i}$
intersects all stable manifolds transversely.
\end{proofof}

Minor modifications of the above argument provide a slightly
stronger result

\begin{theorem}
\label{morse-smale-f} Assume a circle action on a compact manifold
$M$ is semi-free and has isolated fixed points. Then for every
Morse function $f$ there exists a gradient-like vector filed $\xi$
such that $(f,\xi)$ is Morse-Smale.
\end{theorem}

\begin{remark} {\rm As mentioned in the introduction, it is not true that any compact manifold
with a circle action has an invariant Morse-Smale function.
}
\end{remark}

\section{Cancellation Theorems.}
\label{section-cancell}

The following example is an illustration to the cancellation
theorems.

\begin{example}
\label{example} {\rm Let $S^1$ act on $\mathbb R^2=\mathbb C$ with
the weight $1$. Consider a function \mbox{$f(z)=v(|z|)$} where
$v(r)$ is a real valued function which is equal to $r^2$ near
zero, strictly increases as~$r$ increases from $0$ to $1$, equals
to $2r-r^2$ near $1$, and strictly decreases as $r$ goes from~$1$
to infinity. Then $f$ is an invariant Morse function with a
fixed point at the origin of index~$0$ and a critical circle at
the unit circle of index~$1$. An invariant gradient-like vector
field is given by $(x\rho(|z|), y\rho(|z|))$ at the point
$z=x+\sqrt{-1}y$, where $\rho(r)=-1$  near $0$, $\rho(r)<0$ when
$0<r<1$, $\rho>0$ when $r>1$, and $\rho=\frac{1-r}{r}$ near~$1$.

Changing $f$ to $2-|z|^2$ and the gradient-like vector field to
$(x,y)$ cancels the critical circle and increases the index of the
origin to $2$. See Figure~\ref{figure}.} \hfill $\square$
\end{example}

\begin{theorem}
\label{cancell-one} {\rm (The preliminary cancellation theorem)}
Assume $f$ is an invariant Morse function on a compact manifold
$M$ with semi-free circle action and $\xi$ is a gradient-like
vector field on $M$. Assume $a<b$ and $M_{(a,b)}$ contains a
single fixed point $p$ of index~$\lambda$ and a single critical
circle $s$ of index~$\lambda+1$. Let $T=W^+_p\cap W^-_s$ be a
single disc. Then there exists an invariant Morse function $f'$
identical to $f$ outside $M_{(a,b)}$, such that $M_{(a,b)}$ contains
no critical circles of $f'$ and the index of the fixed point $p$ is
$\lambda+2$.
\end{theorem}

\proof  Let $S^1$ act on $\mathbb C^n$ with the weight
$(1,\dots,1)$. Define
$$
h(z_1,\dots,z_n)=f(p)
+(f(s)-f(p))v(|z_1|)-|z_2|^2-\dots-|z_{\lambda+1}|^2
+|z_{\lambda+2}|^2+\dots+|z_n|^2,
$$
where $v$ is the function from Example \ref{example}. Moreover,
let $\eta$ be the vector field given by
$$
(x_1\rho(|z_1|), y_1\rho(|z_1|),
x_2,y_2,\dots,x_{\lambda+1},y_{\lambda+1},
-x_{\lambda+2},-y_{\lambda+2},\dots,-x_n,-y_n)),
$$
where $\rho$ is the function from Example \ref{example}. Denote by
$\tilde p$ the origin of $\mathbb C^n$ and by $\tilde s$ the unit
circle inside $\mathbb C\times 0\times\dots\times 0$.

\smallskip

\begin{claim} \label{claim} \rm{There exists an invariant diffeomorphism $g$ which maps an
invariant neighborhood $V$ of the disc
$D=\{(z_1,0,\dots,0);|z_1|\leq 1\}$ onto an invariant neighborhood
$U$ of~$T$ with $g^*(f)=h$ and $g_*(\eta)=  \xi$.}
\end{claim}

%

\begin{proofof}{Claim \ref{claim}} There exist diffeomorphisms $g_p$
and $g_s$ from neighborhoods $V_{\tilde p}$ of~$\tilde p$ and
$V_{\tilde s}$ of $\tilde s$ onto neighborhoods $U_p$ of $p$ and
$U_s$ of $s$ which take $D$ onto~$T$. Let~$\varepsilon$ be small
enough so that $f(p)+\varepsilon<f(s)-\varepsilon$. Then set $g$
equal to $g_p$ and to $g_s$
outside~$M_{[f(p)+\varepsilon,f(s)-\varepsilon]}$. By following the
trajectories of $\xi$ and $\eta$ the map $g$ extends to~$D$.

Let $L_{\tilde p}=V_{\tilde p, f(p)+\varepsilon}$ (recall our
conventions $V_{\tilde p, f(p)+\varepsilon}= V_{\tilde p}\cap
h^{-1}(f(p)+\varepsilon)$) and \mbox{$L_{\tilde s}=V_{\tilde s,
f(s)-\varepsilon}$.} Set $D_{\tilde p}=D\cap
h^{-1}(f(p)+\varepsilon)$ and
$D_{\tilde s}=D\cap h^{-1}(f(s)-\varepsilon)$. By following
the
trajectories of $\eta$ we can define an invariant diffeomorphism
$j_1$ from a neighborhood~$W_{\tilde s}$ of $D_{\tilde s}$
inside~$L_{\tilde s}$ onto a neighborhood $W_{\tilde p}$
of~$D_{\tilde p}$ inside~$L_{\tilde p}$. At the same time, if we
apply~$g_s$, then follow the trajectories of $\xi$ and then
apply~$g^{-1}_p$ we get another diffeomorphism~$j_2$ from a
neighborhood
$W'_{\tilde s}$ of $D_{\tilde s}$ inside $L_{\tilde s}$ onto a
neighborhood $W'_{\tilde p}$ of~$D_{\tilde p}$ inside~$L_{\tilde
p}$.

The invariant diffeomorphisms $j_1$ and $j_2$ define
diffeomorphisms $\bar j_1$ and $\bar j_2$ of $W_{\tilde s}/S^1$
to~$W_{\tilde p}/S^1$ and of $W'_{\tilde s}/S^1$ to $W'_{\tilde
p}/S^1$. Using \cite[Theorem~5.6]{m} we can construct an isotopy
of identity $\bar \psi_t$ of a neighborhood $W/S^1$ of $D_{\tilde
p}/S^1$ inside $L_{\tilde p}/S^1$ which intertwines~$\bar j_1$ and
$\bar j_2$, that is $\bar \psi_1(\bar j_1)=\bar j_2$. Since we can
assume $W/S^1$ is contractible, there exists an identification
$W=S^1\times W/S^1$ and hence a lift of $\bar \psi_t$ to an
invariant isotopy $\psi_t$ which intertwines $j_1$ and $j_2$ on
$W$. Since $\psi_t$ is identity outside some compact set the
isotopy~$g_p(\psi_t)$ extends to the whole~$M_{f(p)+\varepsilon}$.

Use Lemma \ref{skew} to perturb $\xi$ on
$M_{[f(p)+\varepsilon,f(s)-\varepsilon]}$ to guarantee that
$j_1=j_2$ for the newly constructed gradient-vector field. The
claim is proved, since $g$ extends to an open neighborhood of $D$
by following the trajectories of $\xi$ and~$\eta$.
\end{proofof}
%

The next step in the proof of Theorem~\ref{cancell-one} is to change
$\xi$ to
$\xi'$ inside $U$
to guarantee that there is only one critical point of $\xi'$ inside
$M_{[a,b]}$.

By argument identical to the proof of Assertion 1 in the proof of
\cite[Theorem 5.4]{m}, there exists an invariant neighborhood
$U''$ of~$T$ inside an invariant neighborhood $U'$ of~$T$ such
that $\overline{U'}\subset U$ and no trajectory of $\xi$ enters
$U''$ then exits $U'$ and then returns to~$U''$.

Change $\xi$ to $\xi'$, such that $\xi=\xi'$ outside of $U''$,
while inside $U''$
$$
\xi'=g_*(v_1(x_1,\tau),v_2(y_1,\tau),x_2,y_2,\dots,
x_{\lambda+1},y_{\lambda+1},
-x_{\lambda+2},-y_{\lambda+2},\dots,-x_n,-y_n)),
$$
where $\tau=|z_2|^2+\dots+|z_n|^2$, $v_1(x_1,0)=x_1$,
$v_2(y_1,0)=y_1$. Moreover, $v_1$ and $v_2$ are chosen such that
$\xi'=\xi$ outside a compact neighborhood of $T$ inside $U''$.
Clearly, $\xi'$ has only one critical point, which is $p$.

Let us show that there is no trajectory $\psi(t,x)$ of $\xi'$
through $x$ which lies completely inside~$U'$. Let~$x=g(x_1^\circ,
y_1^\circ,\dots, x_n^\circ, y_n^\circ)$. If one of the coordinates
$x_2^\circ, y_2^\circ,\dots, x_n^\circ, y_n^\circ$ is not zero,
then this coordinate of $\psi(t,x)$ must increase or decrease
exponentially and hence this trajectory eventually exits $U'$
either as $t$ increases or decreases. Otherwise, if only $x_1$ or
$y_1$ are not zero, then again the trajectory will exit $U'$ as
$t$ increases, since
$\xi'(x_1,y_1,0,\dots,0)=c(x_1,y_1,0,\dots,0)$ for some
positive~$c$.

We now show $\xi'$ has no periodic trajectories and any trajectory
is inside~$W^+_p$ or~$W^-_p$, or connects $M_b$ to $M_a$ (at this
point we forget about $\xi$, so that $W^+_p$ and $W^-_p$ are
the stable and unstable manifolds of $p$ of $\xi'$).
Indeed, if a trajectory of $\xi'$ does not intersect~$U''$, it
connects $M_b$ to $M_a$, since it coincides with a trajectory of
$\xi$. If, on the other hand, a trajectory $\psi(t,x)$ of $\xi'$
through a point $x$ inside $U''$ intersects $U'$ and is not inside
$W^-_p$ or~$W^+_p$, then there should be two nonzero coordinates
of $x$ of indices $i,j$ with $i\leq \lambda+1<j$. Such a
trajectory must exit $U'$ as $t$ increases and as $t$ decreases
and never return to $U'$. Hence this trajectory connects $M_a$ to
$M_b$. This finishes the construction of $\xi'$.

It remains to construct an invariant Morse function for which
$\xi'$ is a gradient-like vector field. To do this we will define
an invariant Morse function $f'$ on $M_{[a,b]}$ with only one
critical point $p$ and with $f'(M_a)=a$, $f'(M_b)=b$ and $f(p)=c$
for $a<c<b$. Then set $f'=f$ outside of $M_{[a,b]}$ and, by the
smoothing procedure described in the proof of Theorem \ref{morse},
smoothen $f'$ in a neighborhood of $M_a$ and $M_b$.

Let $V$ be an invariant neighborhood of $p$ inside a canonical
chart of~$\xi'$. Choose $\varepsilon$ with $a<c-\varepsilon$ and
$c+\varepsilon <b$. We can clearly define $f'$ on
$V_{[c-\varepsilon, c+\varepsilon]}$.

Choose invariant neighborhoods $V^+, V^-$  of $W^+_p, W^-_p$ which
are preserved by the flow of $\xi'$. Then $f'$ can be define on
$V^+_{[c+\varepsilon, b]}$ and $V^-_{[a,c-\varepsilon]}$. Using
again the smoothing argument~$f'$ can be defined on a certain
neighborhood~$W$ of $W^+_p\cup W^-_p$. Since there are no critical
points on $M_{[a,b]}-(W^+_p\cup W^-_p)$, there exists an invariant
function $g$ on $M_{[a,b]}-(W^+_p\cup W^-_p)$ with $\xi'(g)<0$,
$g(M_a)=a$ and $g(M_b)=b$.

Let $\mu$ be an invariant function on $M_a$ which is equal to $0$
outside $M_{a}\cap V^-$ and equal to $1$ in a neighborhood of
$M_a\cap W^-_p$. By following the trajectories of $\xi'$ we get a
map $M_{[a,b]}-(W^+_p\cup W^-_p)\to M_a-W^-_p$. Pull back $\mu$ to
a function $\nu$ on $M_{[a,b]}-(W^+_p\cup W^-_p)$, extend it to
$M_{[a,b]}$ by setting $\nu=1$ on $W^+_p\cup W^-_p$. Then the
function $\nu f'+(1-\nu)g$ is the invariant Morse function on
$M_{[a,b]}$ we are looking for.
\endproof

Let $S^1$ act freely on a manifold $V$. Assume two invariant
submanifolds $N$ and $N'$ of dimensions $r+1$ and $s+1$ intersect
transversely and $r+ s=\dim V-1$. Then $N$ and~$N'$ intersect
along some circles $s_1,\dots,s_k$. Assume $N$ and the normal
bundle $\nu(N')$ are oriented. At a point $x_i\in s_i$ choose a
positively oriented frame $\xi_0,\dots,\xi_r$ on $N$ such that
$\xi_0$ is the restriction to $x_i$ of the infinitesimal vector
field associated to the circle action on $s_i$. Since $N$ and $N'$
intersect transversely, $\xi_1,\dots,\xi_r$ form a basis of the
normal bundle to $N'$ at $x_i$. The intersection number of $N$ and
$N'$ at $s_i$ is defined to be $+1$ if~$\xi_1,\dots,\xi_r$ is a
positively oriented basis and $-1$ otherwise. This definition does
not depend on the choice of the point $x_i$ inside of $s_i$. The
intersection number $N\cdot N'$  of $N$ and $N'$ is the sum of all
the intersection numbers at the circles $s_i$.

If $V$ is oriented, then we say that orientations on $N$ and on
$\nu(N)$ are compatible with respect to the orientation on $V$ if
at any point on $N$ a positively oriented frame on $N$ followed by
a positively oriented frame on $\nu(N)$ gives a positively
oriented frame on $V$. If compatible orientation on $N$ and
$\nu(N)$ as well as on $N'$ and $\nu(N')$ are used to
define~$N\cdot N'$ and $N'\cdot N$ then $N\cdot N'=\pm N'\cdot N$
depending on parity of $rs$.

Before we can prove the cancellation theorem we need two
preliminary results:

\begin{theorem}
\label{cancell-two} Let $N$ and $N'$ be smooth closed invariant
transversely intersecting submanifolds of dimensions $r+1$ and
$s+1$ in the $(r+s+1)$-manifold $V$ on which $S^1$ acts freely.
Suppose $N$ and the normal bundle to $N'$ are oriented. Assume
$r+s\geq 5$, $s\geq 3$ and, in case $r=2$, suppose the inclusion
map $\pi_1((V-N')/S^1)\to\pi_1(V/S^1)$ is injective.

Let $s_1,s_2\in N\cap N'$ be circles with opposite intersection
numbers. Assume there exist two invariant embeddings of $S^1\times
[0,1]$ into $N$ and into $N'$ which connect $s_1$ to $s_2$, such
that both embeddings miss $N\cap N'-\{s_1,s_2\}$. The composition
of two embeddings provides a map of $S^1\times S^1$ into $V$.
Assume $S^1\times S^1$ is invariantly contractible to a circle
inside~$V$.

With these assumptions there exists an invariant isotopy of
identity $h_t$ on $V$, such that $h_t$ is identity near $N\cap
N'-\{s_1,s_2\}$ and $h_1(N)\cap N'=N\cap N'-\{s_1,s_2\}$.
\end{theorem}

\proof We deliberately stated this theorem in the form very
similar to \cite[Theorem~6.6]{m}. Actually, \cite[Theorem 6.6]{m},
applied to manifolds $V/S^1$, $N/S^1$, and $N'/S^1$ provides an
isotopy
$\tilde h_t$ of $V/S^1$ which is identity near $N/S^1\cap
N'/S^1-\{s_1/S^1,s_2/S^1\}$ and such that $\tilde h_1(N/S^1)\cap
N'/S^1=N/s^1\cap N'/S^1-\{s_1/S^1,s_2/S^1\}$. So, it is enough to
show that $\tilde h_t$ can be lifted to an invariant isotopy $h_t$
on $V$.

In the proof  of \cite[Theorem 6.6]{m} the constructed  isotopy
$\tilde h_t$ is identity outside some contractible open set $U$ in
$V$. (If $D$ is the disc in $V/S^1$, whose boundary consists of
two arcs connecting $s_1/S^1$ to $s_2/S^1$, then $U$ is a
neighborhood of $D$.) Set $h_t$ to be equal to identity outside of
$\pi^{-1}(U)$ (where $\pi:V\to V/S^1$). Choosing an invariant
trivialization $\pi^{-1}(U)=U\times S^1$ allows us to lift $\tilde
h_t$ to $h_t$ on $U$.
\endproof

\begin{lemma}
\label{fundamental-group} Given  a semi-free circle action  on a
compact
connected manifold $M$, let~$f$ be an invariant Morse
function on $M$ and $\xi$ a gradient-like vector field. For
regular values $a<b$ of $f$,
suppose $M_{[a,b]}$ contains exactly one critical point $p$ and $M_a$
is connected then
$$
\pi_1(M_a/S^1)=\pi_1(M_b/S^1)=\pi_1((M_a-W^-_{p,a})/S^1)=
\pi_1((M_b-W^+_{p,b})/S^1).
$$
\end{lemma}

\proof Since $p$ is of even index, $M_b$, $M_a-W^-_{p,a}$ and
$M_b-W^+_{p,a}$ are connected and it makes sense to talk about
their fundamental groups. The trajectories of $\xi$ provide an
invariant diffeomorphism of $M_a-W^-_{p,a}$ and $M_b-W^+_{p,b}$,
in particular, $\pi_1((M_a-W^-_{p,a})/S^1)=
\pi_1((M_b-W^+_{p,b})/S^1).$

Let $U$ be an invariant tubular neighborhood of $W^-_{p,a}$ inside
$M_a$. Then $\pi_1(U/S_1)= \pi_1(W^-_{p,a}/S^1)=1$, since
$W_{p,a}/S^1$ is homeomorphic to a complex projective space.
Moreover, \mbox{$\pi_1(U-W^-_{p,a}/S^1)=1$}, since
$U-W^-_{p,a}/S^1$ is homotopy equivalent to an even sphere bundle
over a complex projective space. Hence, van Kampen's theorem
applied to the covering of $M_a/S^1$ by $(M_a-W^-_{p,a})/S^1$ and
$U/S^1$ immediately leads to
$\pi_1(M_a/S^1)=\pi_1((M_a-W^-_{p,a})/S^1)$. Analogously, we can
prove $\pi_1(M_b/S^1)=\pi_1((M_b-W^+_{p,b})/S^1)$. \endproof

\begin{proofof}{Theorem \ref{cancell-theorem}} Taking into account that
$\lambda$ is even, there are five different cases:

\noindent{\bf Case 1:} {\underline{$\lambda=0$.}} Then $N'$ is two
copies of a circle and $N\cdot N'=\pm 1$ implies $N\cap N'$ is
just one circle hence Theorem \ref{cancell-one} applies.

\noindent{\bf Case 2:} {\underline{$\lambda=\dim M-2$.}} Then $N$
is a circle and Theorem \ref{cancell-one} applies again.

\noindent{\bf Case 3:} {\underline{$\dim M-4>\lambda\geq 4$, $\dim
M\geq 8$, $M_a/S^1$ is simply-connected, $M_a$ is connected.}} Set
$r+1=\dim N=\dim M-\lambda-1$ and $s+1=\dim N' =\lambda+1$. Then
$r+s=\dim M-2\geq 5$, $s=\lambda\geq 3$, and $r=\dim
M-\lambda-2\geq 3$. Lemma \ref{fundamental-group} implies
$\pi_1(M_c/S^1)=1$, hence Theorem \ref{cancell-two} applies to
this situation. Repeatedly apply  Theorem \ref{cancell-two} to
find an invariant homotopy $h_t$ of $M_c$ to guarantee $h_1(N)\cap
N'$ is a single circle. Then we can apply Lemma~\ref{skew} to find
new $(f',\xi')$ for which $W^+_p$ intersect $W^-_s$ along a single
disc. Finally we apply the preliminary cancellation theorem to
finish this case.

\noindent{\bf Case 4:} {\underline{$\lambda=\dim M-4$, $\dim M\geq
8$, $M_a/S^1$ is simply-connected, $M_a$ is connected.}} Define
$r+1=\dim N=\dim M-\lambda-1 =3$, $s+1=\dim N'=\lambda+1$, so that
$r=2$, $s\geq 3$ while $r+s=\dim M-2\geq 5$. So, to apply Theorem
\ref{cancell-two} and then Theorem~\ref{cancell-one} as in {\bf
Case~3}, it is enough to show $\pi_1((M_c-N')/S^1)=
\pi_1(M_c/S^1)=1$. Notice $\pi_1(M_c/S^1)=1$ by
Lemma~\ref{fundamental-group}. To show $\pi_1((M_c-N')/S^1)= 1$,
identify $N'/S^1$ with the sphere of dimension $\dim M-3$ and then
use an argument analogues to the proof of
Lemma~\ref{fundamental-group}.

\noindent{\bf Case 5:} {\underline{$\lambda=2,\dim M\geq 8$,
$M_a/S^1$ is simply-connected , $M_a$ is connected.}} Here we
interchange the roles of $N$ and $N'$ in
Theorem~\ref{cancell-two}. Namely, set $r+1=\dim  N'=3$ and
$s+1=\dim N =\dim M-3$. Then $r=2,$ $r+s\geq 5$ and $s\geq 3$, so
to apply Theorem~\ref{cancell-two} and the preliminary
cancellation theorem as in {\bf Case~3}, it is enough to check
$\pi_1((M_c-N)/S^1)=\pi_1(M_c/S^1)$, which immediately follows
from Lemma \ref{fundamental-group}.
\end{proofof}

\section{Equivariant Cohomology and Morse Theory.}
\label{cohomology}

This section discusses properties of equivariant cohomology and
its connections to Morse theory. Let $ES^1$ be a contractible
space with a free $S^1$ action. (For example,~$ES^1$ can be the
infinite dimensional sphere $S^\infty$.) The equivariant
cohomology is defined by
\begin{eqnarray*}
H_{S^1}^*(M,\mathbb Z)=H^*((M\times  ES^1)/S^1,\mathbb Z).
\end{eqnarray*}
We only use cohomology with integer coefficients and write
$H_{S^1}^*(M)= H_{S^1}^*(M,\mathbb Z)$.

If a circle action on $M$ is free then $H^*_{S^1}(M)=H^*(M/S^1)$.

The equivariant cohomology of a point with a trivial circle action
is equal to the cohomology of the space $BS^1=ES^1/S^1$. (If
$ES^1=S^\infty$, then $BS^1$ is the infinite dimensional
projective space $\mathbb C P^\infty=S^\infty /S^1$.) It is known
that $H^*_{S^1}(pt)=H^*(BS^1)$ is the ring of polynomials in one
variable $\mathbb Z[u]$ where the generator $u$ has degree $2$.

The map $M\to pt$ produces the pullback map in equivariant
cohomology $H^*_{S^1}(pt)\to H^*_{S^1}(M)$ which equips the
equivariant cohomology ring with a $\mathbb Z[u]$-module
structure.

A natural inclusion of $M$ into $(M\times ES^1)/S^1$ gives a map
from the equivariant cohomology to the regular cohomology $\kappa:
H^*_{S^1}(M)\to H^*(M).$

Given an invariant Morse function $f$ on $M$ and $a<b$, the exact
sequence of the tuple~$(M_{\leq a},M_{\leq b})$ is
\begin{equation}
\label{exact-sequence} ...\to H^*_{S_1}(M_{\leq b}) \to
H^*_{S_1}(M_{\leq a}) \stackrel{\delta}{\to} H^{*+1}_{S_1}(M_{\leq
b}, M_{\leq a}) \to H^{*+1}_{S_1}(M_{\leq b})\to ...
\end{equation}
To understand this exact sequence (especially the map $\delta$)
better, let us describe the cohomology $H^*_{S_1}(M_{\leq
b},M_{\leq a})$ in the case when $M_{(a,b)}$ contains a single
critical point $p$ of index $\lambda$. By the classical results in
invariant Morse theory (see \cite{ab}) we know that
$$
H^*_{S_1}(M_{\leq b},M_{\leq a})=H^*_{S^1}(W^-_{p, \geq a}
W^-_{p,a}) = H^*_{S^1}(D^\lambda, S^{\lambda-1}),
$$
where $D^\lambda$ is the unit disc inside $\mathbb
C^{\frac{\lambda}{2}}$ on which $S^1$ acts with the weight
$(1,\dots,1)$ and $S^{\lambda-1}=\partial D^{\lambda}$. By Thom
isomorphism $H^*_{S^1}(D^\lambda,
S^{\lambda-1})=H^{*-\lambda}_{S^1}(pt)$. Choosing a
generator~$\tau_p$ of $H^{\lambda}_{S^1}(D^\lambda,
S^{\lambda-1})$ we get a $\mathbb Z[u]$-module identification
$$
H^*_{S_1}(M_{\leq b},M_{\leq a})=\tau_pH^*_{S^1}(pt).
$$
We call $\tau_p$ \emph{a Thom class of $p$}. Notice that choosing
$\tau_p$ is equivalent to picking an orientation on $W_p^-$.

Similarly, if $M_{(a,b)}$ has a single critical circle $s$ of
index $\lambda$ then
$$
H^*_{S^1}(M_{\leq b},M_{\leq a})=\tau_s H^*_{S^1}(S^1)=\tau_s
H^*(pt)
$$
and $\tau_s\in H^\lambda_{S^1}(M_{\leq b},M_{\leq a})$ is \emph{a
Thom class of $s$.}

\begin{lemma}
\label{restatement} Assume the action of $S^1$ on a compact
oriented manifold $M$ is semi-free and $(f,\xi)$ is an invariant
Morse-Smale function on $M$. Assume $a<b$ and $M_{(a,b)}$ contains
a critical point $p$ of index $\lambda$ and a critical circle $s$
of index $\lambda+1$. Let \mbox{$f(p)<c<f(s)$}. Suppose $M_c$ is
oriented and $N=W^+_{p,c}$, $\nu(N)$ as well as $N'=W^-_{s,c}$,
$\nu(N')$ are equipped with compatible orientations.

Consider the boundary map of the exact sequence of the triple
$(M_{\leq b}, M_{\leq c}, M_{\leq a})$
$$
\delta:H^*_{S^1}(M_{\leq c}, M_{\leq a})\to H^{*+1}_{S^1}(M_{\leq
b}, M_{\leq c}).
$$
This map takes \mbox{$H^{\lambda}_{S^1}(M_{\leq c}, M_{\leq
a})\cong \tau_p\mathbb Z$} to $H^{\lambda+1}_{S^1}(M_{\leq b},
M_{\leq c})\cong \tau_s\mathbb Z$. If $\delta(\tau_{p}) =c \tau_s$
then~$N \cdot N'=\pm c$.
\end{lemma}

\proof Use excision to identify
$$
H^*_{S^1}(M_{\leq c}, M_{\leq a})=H^*_{S^1}(M_{[a,c]}, M_a)
$$
and
$$
H^*_{S^1}(M_{\leq b}, M_{\leq c})=H^*_{S^1}(M_{[c,b]}, M_{ c}).
$$

The map
$$
\kappa:H^*_{S^1}(M_{[a,c]}, M_a)\to
H^*(M_{[a,c]}, M_a)
$$
is an isomorphism in degree $\lambda$. Moreover, by Poincare
duality
$$
H^\lambda(M_{[a,c]}, M_a)=H_{\dim M-\lambda}(M_{[a,c]},M_c).
$$

Since $S^1$ acts freely on $M_{[c,b]}$ we have
$$
H^*_{S^1}(M_{[c,b]}, M_{ c})=H^*(M_{[c,b]}/S^1, M_{ c}/S^1).
$$
By Poincare duality
$$
H^{\lambda+1}(M_{[c,b]}/S^1, M_{ c}/S^1)\cong H_{\dim M
-\lambda-2}(M_{[c,b]}/S^1, M_b/S^1).
$$

Thus, we view $\delta$ as a map
$$
\delta:H_{\dim M-\lambda}(M_{[a,c]},M_c) \to H_{\dim M
-\lambda-2}(M_{[c,b]}/S^1, M_b/S^1).
$$
This map has a geometric interpretation. The generator of $H_{\dim
M-\lambda}(M_{[a,c]},M_c)$ is given by the cycle $W^+_{p,[a,c]}$,
then $\delta$ sends it to the cycle $(\partial
W^+_{p,[a,c]})/S^1=W^+_{p,c}/S^1=N/S^1$ inside the tuple
$(M_{[c,b]}/S^1, M_b/S^1)$.

\cite[Lemma~7.2]{m} used in our situation states that the homology
class of $N/S^1$ is a~generator of $H_{\dim
M-\lambda-2}(M_{[c,b]}/S^1, M_b/S^1)$  multiplied by $N\cdot N'$
which proves the~lemma.
\endproof

\begin{theorem}
\label{basis} Assume the action of $S^1$ on a compact oriented
manifold $M$ is semi-free and $(f,\xi)$ is an invariant
Morse-Smale function. For $a<b$ suppose $M_{[a,b]}$ is connected.
Assume the set of critical points in $M_{(a,b)}$ consists of
critical points $p_1,\dots,p_k$ of the same index $\lambda$ with
\mbox{$2\leq\lambda\leq \dim M-2$} (or critical circles
$s_1,\dots,s_k$ of the same index~$\lambda$ with
\mbox{$3\leq\lambda\leq \dim M-2$}). Given a basis
$\tau_1,\dots,\tau_k$ of $H^\lambda_{S_1}(M_{\leq b},M_{\leq a})$
there exists an invariant Morse-Smale function $(f',\xi')$ which
coincides with $(f,\xi)$ outside $M_{[a,b]}$, has the same set of
critical points as $(f,\xi)$ and after choosing the proper
orientations on $W^-_{p_i}$ the Thom class of each $p_i$ is
$\tau_i$ (or the Thom class of each $s_i$ is $\tau_i$).
\end{theorem}

\proof Consider the case when $M_{(a,b)}$ contains only critical
points $p_1,\dots,p_k$. Forget for a moment about the circle
acting on our manifold and try to prove the theorem for the
regular cohomology groups. Then using excision and Poincare
duality we can identify
$$
H^\lambda(M_{\leq b},M_{\leq a})=H_{n-\lambda}(M_{[a,b]}, M_b).
$$
Moreover, the Thom classes of the fixed point $p_i$ are given by
the cycles $W^+_{p_i,[a,b]}$. Translating our theorem into a
statement about a basis in homology leads to a statement identical
to  . In particular, there exists a Morse-Smale function $(\tilde
f, \tilde \xi)$ on $M$ which coincides with $(f,\xi)$ outside
$M_{[a,b]}$ and the Thom classes (in regular cohomology) of
$p_i$'s are exactly $\kappa(\tau_i)$'s

Notice that the map $\kappa$ identifies $H^\lambda_{S^1}(M_{\leq
b},M_{\leq a})$ with $H^\lambda(M_{\leq b},M_{\leq a})$. So, if we
show $(f',\xi')$ can be chosen invariantly with respect to the
circle action, the theorem is proved.

Without loss of generality we may assume $(\tilde f,\tilde \xi)$
coincides with $(f,\xi)$ in a neighborhood of every fixed point,
this follows from an examination of the proof of
\cite[Theorem~7.6]{m}.

Let us average $(\tilde f, \tilde \xi)$ with respect to the circle
action to produce invariant function $f'$ and invariant vector
field $\xi'$. Then $(f',\xi')$ is identical to $(f,\xi)$ outside
of $M_{[a,b]}$ and in neighborhoods of the fixed points
$\{p_1,\dots, p_k\}$. Moreover, it is obvious that $\xi'(f')<0$ on
$M-M^{S^1}$. Hence $f'$ is an invariant Morse function and $\xi'$
is its invariant gradient-like vector field. Moreover
$\kappa(\tau_{p_i})$ are the same as for $(\tilde f, \tilde \xi)$.
Hence the Thom classes $\tau_{p_i}$ form the required basis.

Finally, to make $(f',\xi')$ into an invariant Morse-Smale
function we need to apply Theorem \ref{morse-smale} and perturb
$\xi'$ without changing the Thom classes $\tau_{p_i}$.

Let us now outline the proof in the case when $M_{(a,b)}$ contains
only critical circle. Since the circle acts freely on $M_{[a,b]}$,
we can make an analogous statement on $M_{[a,b]}/S^1$ and then
using Poincare duality to translate it into a statement in
homology. This statement will follow from \cite[Theorem~7.6]{m},
which is proved by constructing a series of isotopies of identity
on $M_a$ and then changing the vector field according to these
isotopies. Inspection of the proof of \cite[Theorem~7.6]{m} shows
that each of these isotopies is identity outside a contractible
set and hence can be lifted to $M_{[a,b]}$. This adopts the proof
of \cite[Theorem~7.6]{m} to give a proof of our theorem.
\endproof

\section{Proof of the main theorem.}

\begin{lemma}
\label{formality} Assume a circle action on a compact manifold $M$
is torsion-free and has isolated fixed points $p_1,\cdots, p_m$.
Choose invariant neighborhoods $U_i$ of $p_i$ diffeomorphic to
unit balls and set $U=\cup_iU_i$. Then the homology groups
$H_*((M-U)/S^1,U/S^{1})$ are $\ZZ$-torsion-free.
\end{lemma}

\proof By Poincare duality
$$
H_*((M-U)/S^1,\partial U/S^1)=H^{\dim M-1-*}((M-U)/S^1).
$$
Since
$$
H^{*}((M-U)/S^1)=H^{*}_{S^1}(M- U)=H^{*}_{S^1}(M-M^{S^1}),
$$
it is enough to show that the cohomology groups
$H^{*}_{S^1}(M-M^{S^1})$ are torsion-free.

Consider the Mayer-Vietoris sequence for $M=U\cup(M-M^{S^1})$:
$$
... \to H^*_{S^1}(M)\to H^*_{S^1}(U)\oplus H^*_{S^1}(M-M^{S^1})\to
H^*_{S^1}(U\cap (M-M^{S^1}))\to ...
$$
An easy computation in equivariant cohomology implies the
surjectivity of the map
$$
\oplus_{i}H^*_{S^1}(p_i) \cong H^*_{S^1}(U) \to H^*_{S^1}(U\cap
(M-M^{S^1}))\cong H^*(\partial U/S^1).
$$
Hence the Mayer-Vietoris  sequence breaks up into short exact
sequences. Since $H^*_{S^1}(M)$ and $H^*_{S^1}(U\cap (M-M^{S^1}))$
are torsion-free, the same holds for
\mbox{$H^{*}_{S^1}(M-M^{S^1})$.}
\endproof

An invariant Morse function $f$ is \emph{self-indexing} if
$f(q)=\sigma(q)$ for every $q\in Crit_c(f)$, in particular, the
critical values of $f$ are between zero and $\dim M$. For~an
integer~$k$, define $M_k=M_{\leq k+\frac{1}{2}}$. Then
$H_{S^1}^*(M_k, M_{k-1})$ is generated as a $\ZZ[u]$ module by the
Thom classes of critical points and circles of index $k$.  Then
$f$ is perfect on~$M_k$ if~\mbox{$\mathcal P(M_k)=\mathcal
M(M_k,f)$}. The exact sequence of the triple $(M_{k+1}, M_k,
M_{k-1})$ defines the~map
\begin{equation*}
\delta_k : H^*_{S^1}(M_{k-1}, M_{k-2})\to H^{*+1}_{S^1}(M_{k},
M_{k-1})
\end{equation*}

\begin{proofof}{Theorem \ref{perfect-morse}}
The proof will use induction on $k$ to construct a self-indexing
invariant Morse-Smale function $(f, \xi)$ on $M$, such that
\begin{enumerate}
\item $f$ is perfect on $M_k$,
\item if $s$ is a critical circle of index $k-1$ then
$\delta_k(\tau_s)=0$.
\end{enumerate}

\noindent{\bf Basis of induction for $k=0$.}
Use the notations from Lemma \ref{formality}. Since $M$ is
simply-connected, it is easy to see using van Kampen theorem that
$(M-U)/S^1$ is simply-connected as well. This together with Lemma
\ref{formality} allow to use \cite[Theorem~6.5]{s} which implies
that there exists a perfect Morse function
$$
\tilde f:((M-U)/S^1, (\partial U)/S^1) \to ([1,\infty),1).
$$
Lift $\tilde f$ to $f$ on $M-U$. As in the proof of
Theorem~\ref{morse} extend $f$ to $M$ by setting
$f=\phi_i^*(|z|^2)$ on $U_i$ and smoothing it along $\partial U$.

By Theorem \ref{morse-smale-f} there exists an invariant
gradine-like vector field $\xi$, such that $(f,\xi)$ is an
invariant Morse-Smale function. By Theorem \ref{order} we can
assume $f$ is self-indexing. Conditions (1) and (2) are satisfied
by $f$, since $\tilde f$ is a perfect Morse function.

\smallskip





\noindent{\bf Induction step for $k=1$.} Assume $f$ satisfies
conditions (1) and (2) for $k=0$.
Let~$p$ be a fixed point of index $0$ such that
$\delta_1(\tau_p)\neq 0$. Then it is easy to see that there exists
a critical circle $s$ of index $1$ such that the closure of the
intersection $W^-_s\cap W^+_p$ is a single disc. By Theorem
\ref{cancell-one} we can cancel the critical circle $s$ by
increasing the index of $p$ by $2$. We can continue this process
until the map $\delta_1$ is nontrivial. After this $(f,\xi)$ will
satisfy~(1) and (2) for $k=1$.
\smallskip

\noindent{\bf Induction step for even $k\geq 2$.} If $k$ is even
and $(f, \xi)$ satisfies conditions (1) and (2) for~$k-1$, then
$(f,\xi)$ satisfies~(1) and (2) for $k$ as well. Indeed, the
second condition does not depend on~$k$. At the same time, the map
$\delta_k$ is trivial, by condition (2) and the fact that
$M_{k-1}-M_{k-2}$ contains no fixed points. This implies that the
exact sequence of the pair $(M_k, M_{k-1})$ splits into short
exact sequences proving the perfectness of $f$ on $M_k$.

\smallskip

\noindent{\bf Induction step for odd $k\geq 1$.} Assume $(f, \xi)$
satisfies~(1) and (2) for $k-1$. Since the map $\delta_k$ vanishes
on all the Thom classes of critical circles we can assume without
loss of generality that there are no critical circles of index
$k-1$ Assume $p_{i_1},\dots, p_{i_\ell}$ are all the fixed point
of index $k-1$ of $f$.




Choose a basis $\tau_1,\dots, \tau_q$ of the kernel of the map
$\delta_k$ restricted to $H^{k-1}(M_{k-1}, M_{k-2})$. Complete it
to a basis of $H^{k-1}(M_{k-1}, M_{k-2})$ by adding
$\tau_{q+1},\dots, \tau_\ell$. By Theorem~\ref{basis} we can
perturb $(f,\xi)$ on $M_{k-1}-M_{k-2}$ such that
$\tau_{p_i}=\tau_i$ after the perturbation.

Assume $s_1,\dots, s_r$ are the critical circles of $f$ of index
$k$, then their Thom classes form a basis of $H^k(M_k, M_{k-1})$.
The elements $\tilde \tau_{j}=\delta_k(\tau_{q+j})$ of $H^k(M_k,
M_{k-1})$ span the image of $\delta_k$. Complete $\tilde
\tau_{1},\dots\tilde\tau_{\ell-q}$ to a basis of $H^k(M_k,
M_{k-1})$ by adding $\tilde \tau_{\ell-q+1},\dots,\tilde \tau_r$.
By Theorem~\ref{basis} we can alter $(f,\xi)$  on $M_{k}-M_{k-1}$
to guarantee $\tau_{s_i}=\tau_i$. Apply Lemma~\ref{restatement}
together with Theorem~\ref{cancell-theorem} to cancel the critical
circles $s_1,\dots s_{\ell-q}$ by increasing the index of each
$p_{q+i}$ by $2$. This produces a new invariant Morse-Smale
function which satisfies~(1) and (2) for $k$.




The above argument works as long as the assumptions of
Theorems~\ref{cancell-theorem} are satisfied. Checking these
assumptions boils down to making sure that if $f$ is  perfect
on~$M_k$ for $k>1$, then  $(f^{-1}(a))/S^1$ for any regular
$a<k+\frac{1}{2}$ is simply-connected.
It is enough to prove that if the interval $[a,b]$ contains
exactly one critical value then
$\pi_1(f^{-1}(a)/S^1)=\pi_1(f^{-1}(b)/S^1)$.
This follows from Theorem~\ref{order},
Lemma~\ref{fundamental-group} and an analogue of
Lemma~\ref{fundamental-group} for critical circle instead of fixed
point with almost identical proof.
\end{proofof}

\begin{remark} {\rm The dimension restriction, which appears in this
theorem comes from the dimension restrictions in the classical
Morse theory. In particular, if $\dim M$ is either~$4$ or $6$ then
\cite[Theorem~6.5]{s} used in the basis of induction fails.
}
\end{remark}

\end{document}